\documentclass[12pt]{article}
\usepackage[utf8]{inputenc}	
\usepackage[T2A]{fontenc}
\usepackage[english,russian]{babel}
\usepackage{amsmath,amsthm,amsfonts,amssymb}

\newcommand{\ep}{\varepsilon}

\newcommand{\RNumb}[1]{\uppercase\expandafter{\romannumeral #1\relax}} 
\renewcommand{\Re}{\operatorname{Re}} 

\def\le{\leqslant}
\def\ge{\geqslant}

\begin{document}	

\newtheorem{theorem}{Теорема}
\newtheorem*{theorem*}{Теорема}
\newtheorem{lemma}{Лемма}
\newtheorem{note}{Замечание}
\newtheorem{conclusion}{Следствие}
\newtheorem{definition}{Определение}
\renewcommand{\proofname}{Доказательство}
\setlength{\jot}{10pt}

\begin{center}
{{\large {\bf Спектральные асимптотики решений $2\times 2$ системы обыкновенных дифференциальных  уравнений первого порядка}}}
\end{center}
\vskip .5cm
\begin{center}
А.\,П.~Косарев, А.\,А.~Шкаликов\footnote {Исследование поддержано
Российским фондом фундаментальных исследований, грант No~20-11-20261.}
\end{center}

{\bf Ключевые слова}: Асимптотики решений систем обыкновенных дифференциальных уравнений, регулярные краевые задачи, полнота собственных функций.

\bigskip
\bigskip
\section{Введение}
Основная цель этой статьи  --- получить асимптотические представления фундаментальной матрицы решений $2\times 2$ системы дифференциальных уравнений вида
\begin{equation}\label{MainVectorEquation}
	\begin{pmatrix}
		y_1 \\
		y_2
	\end{pmatrix}^{\prime} -
	\begin{pmatrix}
		b_{11} & b_{12} \\
		b_{21} & b_{22}
	\end{pmatrix}
	\begin{pmatrix}
		y_1 \\
		y_2
	\end{pmatrix}
	= \lambda
	\begin{pmatrix}
		a_1 & 0\\
		0 & a_2
	\end{pmatrix}
	\begin{pmatrix}
		y_1 \\
		y_2
	\end{pmatrix}, \,\ \,\ \,\ x \in [0, 1],
\end{equation}
где $\lambda$ - спектральный параметр.

Асимптотические представления будут получены при больших значениях спектрального параметра равномерно по $x\in [0,1]$  в полуплоскостях $\Re\,\lambda \geqslant -\kappa$ и
$\Re\,\lambda \leqslant \kappa$,  где $\kappa$ произвольное вещественное число. В частности, в качестве полуплоскостей могут выступать  левая  и правая полуплоскости $\mathbb C^+$  и $\mathbb C^-$.
Всюду далее будем предполагать, что при некотором $\varepsilon >0$
 \begin{equation} \label{a1}
 a_1 = a_1(x) > 0, \ \  a_2 = a_2(x) < 0 \text{ и } a_1(x) - a_2(x) \geqslant \varepsilon
 \end{equation}
при некотором $\varepsilon >0$.
 В работе будут рассмотрены два варианта гладкости коэффициентов. В первом случае предполагаем
  \begin{equation} \label{a2}
 a_i, \ b_{ij} \in L_1[0,1].
 \end{equation}
Во втором случае предполагаем, что при некотором фиксированном $n\geqslant 1$
 \begin{equation} \label{a3}
 a_1,\ \, a_2,\ \, b_{21},\ \, b_{12} \in W^n_1[0,1], \quad b_{11},\ \, b_{22} \in W^{n-1}_1[0,1],
 \end{equation}
 где $W^k_1[0,1]$ --- соболевское пространство функций $y$,  для которых производные $y^{(s)}$  абсолютно непрерывны при $s\leqslant k-1$,  а $y^{(k)} \in L_1[0,1]$.
При этом предполагаем, что функции $b_{ij}$  могут принимать комплексные значения.

\textit{Матрицей фундаментальной системы решений, или фундаментальной матрицей,} уравнения \eqref{MainVectorEquation} мы называем матрицу
$$
Y(x, \lambda)
=
\begin{pmatrix}
	y_{11}(x, \lambda) & y_{12}(x, \lambda) \\
	y_{21}(x, \lambda) & y_{22}(x, \lambda)
\end{pmatrix},
$$
определенную при $x \in [0, 1]$,  $\lambda \in \mathbb{C}$, $|\lambda| > \lambda_0$, столбцы которой являются независимыми решениями системы \eqref{MainVectorEquation}. Таким образом, матрица $Y(x, \lambda)$ удовлетворяет матричному дифференциальному уравнению
\begin{equation}\label{MainMatrixEquation}
	Y^{\prime}(x, \lambda) = \{\lambda A(x)+B(x)\}Y(x, \lambda)
\end{equation}
при каждом фиксированном $\lambda$. Здесь
 \begin{equation} 
B= \begin{pmatrix}
		b_{11} & b_{12} \\
		b_{21} & b_{22}
	\end{pmatrix}, \quad \	
	A = \begin{pmatrix}
		a_1 & 0\\
		0 & a_2
	\end{pmatrix}.
 \end{equation}

Важную роль в дальнейшем изложении играют матрицы
\begin{equation}\label{ME}
	M(x) =
	\begin{pmatrix}
		e^{\int_0^xb_{11}(t)dt} & 0 \\[0.2cm]
		0 & e^{\int_0^xb_{22}(t)dt}
	\end{pmatrix}, \,\ \,\ \,\
	E(x, \lambda)  = 	\begin{pmatrix}
		e^{\lambda \int_0^x a_{1}(t)dt} & 0 \\[0.2cm]
		0 & e^{\lambda \int_0^x a_{2}(t)dt}
	\end{pmatrix}
\end{equation}

Основной результат работы --- доказать, что существует фундаментальная матрица решений системы \eqref{MainVectorEquation}, которая имеет представление
\begin{equation}\label{R}
	Y(x, \lambda) = M(x)\left(I +  \frac{R^1(x)}{\lambda} + \dots + \frac{R^n(x)}{\lambda^n} + o(1)\lambda^{-n}\right)E(x, \lambda),
\end{equation}
в полуплоскостях $\Pi_{\kappa}^{-} = \{\lambda \in \mathbb{C} | \Re\lambda < \kappa\}$ и $\Pi_{\kappa}^{+} = \{\lambda \in \mathbb{C} | \Re\lambda > -\kappa\}$ для произвольно фиксированного $\kappa\in\mathbb R$. При выполнении условия \eqref{a2} (это соответствует случаю $n=0$) выражение в скобках в  представлении \eqref{R} принимает вид  $(I + o(1))$,  где $o(1) \to0$ равномерно по $x$  при $|\lambda| \to\infty$  в выбранной полуплоскости $\Pi_{\kappa}^{-}$ или  $\Pi_{\kappa}^{+}$. Наиболее важный  момент состоит в том, что при $n\geqslant 1$ мы выписываем явные формулы для матриц $R^m(x),\
1\leqslant m \leqslant n$, участвующих в \eqref{R}. Ранее это не делалось, несмотря на долгую историю исследований на эту тему.

Первые существенные результаты по существованию  асимптотических представлений для фундаментальных решений в секторах комплексной плоскости не только для $2\times 2$  систем, но и для $n\times n$  систем, были получены Я.Д.Тамаркиным  в книге \cite{T}.
При этом предполагалось, что $A= \text{diag}\, (a_1, a_2, \dots, a_n)$, где $a_j$ различные  отличные от нуля числа, а коэффициенты матрицы $B$ имеют  $n+2 $  непрерывных производных. Явные выражения для матриц $R^m$  не выписывались. Другим методом похожие результаты получили Дж.Биркгоф и Р.Лангер \cite{BL}.  В последующие годы на эту тему появилось большое число работ, достаточно подробные  ссылки имеются в недавней работе
Савчука и Шкаликова \cite {SS20}.  В частности, из результатов \cite{SS20}  следует представление \eqref{R},  но только для $n=0$ и при дополнительном условии коллинеарности функций $a_1(x) = -\beta a_2(x), \ \, \beta = \text{const} >0$.

Полученный в этой работе результат является новым даже для системы Дирака (эта система  соответствует случаю $a_1 = -a_2 =1)$, которая исследовалась во многих работах.  Новизна результатов этой работы в сравнении с известными результатами для системы Дирака  состоит в следующем.  Во-первых, вместо чисел $a_1 = -a_2 =1$ в матрице $A$  могут участвовать произвольные функции разных знаков.  Во вторых,  при $n\geqslant 1$  не рассматривался вопрос о минимальных условиях на гладкость коэффициентов матрицы $B$.
 Коэффициенты матрицы $B$ предполагались достаточно гладкими (анализ прежних доказательств позволяет заключить, что для получения представления \eqref{R}  нужно требовать принадлежность коэффициентов матрицы $B$ по меньшей мере пространствам $W^{n+1}_1[0,1]$). В третьих (и этот момент наиболее важен), ранее не выписывались явные формулы для матриц-коэффициентов в  представления \eqref{R}. Эта задача оказалась сложной. Аргументы, свидетельствующие о сложности этой задачи можно найти в книге М.А.Наймарка
 \cite[Гл.1]{Na}.

 Изучение системы  \eqref{MainVectorEquation}  в случае, когда элементы матрицы $A$
 есть функции (вместо констант), не является  простым обобщением, оно продиктовано наличием конкретных задач,  приводящим к необходимости изучать системы с условиями \eqref{a1}. Действительно, система телеграфных  уравнений (см. \cite{TS}, \cite{KS21})  сводится к системе  \eqref{MainVectorEquation}  с коллинеарными функциями $a_1(x) = - a_2(x) >0$ . Нам известны также задачи, которые описывают процессы в металлургии \cite{M}, которые  приводят  к системам вида  \eqref{MainVectorEquation}
 с неколлинеарными функциями $a_1, a_2$  разных знаков.

 В книге \cite{M}   спектральные  задачи для системы  \eqref{MainVectorEquation}
 рассматриваются с краевыми условиями
 \begin{equation}\label{BV}
 y_1(0) = 0, \ \  y_2(1) = 0.
 \end{equation}
 Для  краевых задач, связанных с уравнением   \eqref{MainVectorEquation},
 также, как для системы Дирака, можно определить понятие регулярности. Это понятие можно ввести по-разному, причем для общих $n\times n$  систем. Понятие регулярности связано с задачей о возможности разложения произвольных функций  по собственным функциям краевых задач и  ведет начало от работ \cite{B},
 \cite{T} и \cite{BL}. В более короткой и понятной форме  определение регулярности дано в
 в заметках \cite{S21}, \cite{Sh21}.   В частном случае для системы Дирака
 определение регулярных краевых задач можно найти  в работах \cite{Mar}, \cite{DM}, \cite{SS14}. В \cite{DM} авторы доказали что корневые функции регулярных задач
 для системы Дирака  образуют безусловный базис в пространстве $L_2[0,1] \times L_2[0,1]$   при условии принадлежности коэффициентов матрицы $B$  пространству
 $L_2[0,1]$.  В работе \cite{SS14} это утверждение было доказано при более  общем условии  только суммируемости коэффициентов матрицы $B$,  а в заметке \cite{Sh21}  теорема о безусловной базисности была анонсирована для общих регулярных $n\times n$  систем, названных в  \cite{Sh21}  гиперболическими.

  Из определения регулярности следует, что несмотря на
 простой и естественный вид краевых условий \eqref{BV},  спектральная задача
  \eqref{MainVectorEquation}, \eqref{BV}  регулярной не является. В этом случае
   система корневых функций не будет образовывать базис. В заметке \cite{KS21}
   поясняется, что спектральные свойства задачи  \eqref{MainVectorEquation}, \eqref{BV}  похожи на свойства известной задачи  Редже, которая была изучена во многих статьях (см. \cite{S2} и имеющиеся в этой работе ссылки). Однако задача \eqref{MainVectorEquation}, \eqref{BV}  сложнее, нежели задача Редже. В \cite{KS21} найдено только достаточное условие для полноты системы корневых векторов вместо критерия, найденного в \cite{S2} для задачи Редже.

Важность представлений \eqref{R}  с остатком  $o(1)\lambda^{-n}$  при $n\geqslant 1$  проясняется как раз при исследовании нерегулярных задач. Это продемонстрировано в заметке \cite{KS21}. А именно,  если задача нерегулярна, то оценку ее функции Грина можно провести только с учетом анализа значений матриц $R^m(x)$  в точках $x=0$  и $x=1$.  частности, в \cite {KS21}  показано, что при выполнении  условий
$$
b_{12} \ne 0 \ \, \text{и} \ \, b_{21} \ne 0
$$
 функция Грина задачи \eqref{MainVectorEquation}, \eqref{BV} допускает оценку
 $$
 |G(x, \xi, \lambda)|  \leqslant \text{const} (1+|\lambda|)
 $$
  во всей комплексной $\lambda$-плоскости, но вне кружков фиксированного радиуса с центрами в собственных значениях задачи.  С помощью этой оценки стандартными методами получается доказательство полноты системы корневых функций задачи \eqref{MainVectorEquation}, \eqref{BV} в пространствах $L_p[0,1] \times L_p[0,1]$  при $p\geqslant 1$. Важность многочленных асимптотических представлений можно также оценить  при чтении работы \cite{S83},  где такие представления существенно используются для разных целей.

\bigskip

\section{Асимптотика фундаментальной системы решений}
Далее полезно ввести следующие обозначения
\begin{equation}\label{ab}
\begin{gathered}
	a(x) = a_1(x) - a_2(x), \,\ b(x) = e^{\int_0^x b_{11} (t)- b_{22}(t)dt}, \\ A_i(x) = \int_0^x a_i(t)dt, \,\ \rho(x) = A_1(x) - A_2(x).
\end{gathered}
\end{equation}
Начнём с поиска главного члена в асимптотическом представлении  $Y(x, \lambda)$.
Для этого представим матрицу $B(x)$ в виде суммы диагональной и косодиагольной матриц $B(x)=D(x)+W(x)$,
$$
D(x)
=
\begin{pmatrix}
	b_{11}(x) & 0 \\
	0 & b_{22}(x)
\end{pmatrix},
\,\
W(x)
=
\begin{pmatrix}
	0 & b_{12}(x) \\
	b_{21}(x) & 0
\end{pmatrix}.
$$
Отметим, что матрица $Y_0(x, \lambda) = M(x)E(x, \lambda)$ является решением уравнения \eqref{MainMatrixEquation} в случае, когда $W(x) \equiv 0$ (то есть матрица $B(x)$ - диагональная).  Вид этого решения подсказывает замену, которая значительно упростит наши рассуждения.

Осуществим замену переменных $Y(x, \lambda) = M(x)Z(x, \lambda)E(x, \lambda)$ с неизвестной матрицей $Z(x, \lambda)$. Подставляя $Y(x, \lambda)$ в \eqref{MainMatrixEquation} с учётом равенств $E^\prime = \lambda A E$ и $M^\prime = DM$, получаем
\begin{equation*}
	M^\prime ZE + MZ^\prime E + MZE^\prime = \lambda A MZE+DMZE+WMZE
\end{equation*}
\begin{equation*}
	DMZE + MZ^\prime E + MZ\lambda A E = \lambda A MZE+DMZE+WMZE.
\end{equation*}
\begin{equation*}
	MZ^\prime E= \lambda (AMZ - MZA)E+WMZE.
\end{equation*}

Умножим на $M^{-1}$ слева и на $E^{-1}$ справа и учтём, что диагональные матрицы $M$ и $R$ коммутируют. Тогда
\begin{equation}\label{MainMatrixEquationChange}
	Z^\prime = \lambda(AZ-ZA)+QZ,
\end{equation}
где
\begin{equation}\label{Q}
\begin{gathered}
	Q(x) = M^{-1}(x)W(x)M(x), \,\ q_{11}(x) = q_{22}(x) =0, \\ 
q_{21}(x) = b_{21}(x)b(x), \,\ q_{12}(x) = b_{12}(x)b^{-1}(x).
\end{gathered}
\end{equation}
Будем искать решение \eqref{MainMatrixEquationChange} в виде
\begin{equation}\label{FormalExpansion}
	Z(x, \lambda) = P_0(x)+\frac{1}{\lambda}P_1(x)+\frac{1}{\lambda^2}P_2(x)+\dots
\end{equation}
Приравнивая коэффициенты при степенях $\lambda$ в \eqref{MainMatrixEquationChange}, последовательно определим матрицы $P_0(x), P_1(x)$ и т.д.  Заметим, что далее мы будем строить какую-то фундаментальную матрицу решений. Другой выбор констант интегрирования для элементов $i$-го столбца матрицы $Y(x, \lambda)$ эквивалентен умножению $i$-го столбца матрицы $Y(x, \lambda)$ на некоторый сходящийся ряд $\sum\limits_{s=0}^\infty c_s \lambda^{-s}$ с постоянными коэффициентами $c_s$.

Формальное равенство при $\lambda^1$ дает соотношение
$$
P_0A - AP_0 = 0,
$$
следовательно, $p^0_{12} = p^0_{12} = 0$.

Запишем равенство коэффициентов при $\lambda^0$
$$
P_1A-AP_1+P^\prime_0 = QP_0,
$$
или в покоординатном виде
\begin{equation}\label{lambda_0}
	\left(p_{11}^0\right)^\prime = q_{12}p_{21}^0, \ \ \ \left(p_{22}^0\right)^\prime = q_{21}p_{12}^0, \ \ \ p_{21}^1 = \frac{-\left(p_{21}^0\right)^\prime + q_{21}p_{11}^0}{a}, \ \ \ p_{12}^1 = \frac{-\left(p_{12}^0\right)^\prime + q_{12}p_{22}^0}{-a}.
\end{equation}
Из этих равенств получаем $\left(p_{11}^{0}\right)^\prime = \left(p_{22}^{0}\right)^\prime = 0$. Выбор констант интегрирования в нашей власти, мы положим $p_{ii}^{(0)} = 1$. Тогда $P_0(x) = I$.
Подставляя $p_{11}^{0} = p_{22}^{0} = 1$ в последние два  равенства в \eqref{lambda_0}, находим
$$
p_{21}^1 = \frac{q_{21}}{a} =  \frac{b_{21}b}{a}, \,\ \,\ \,\ p_{12}^1 = \frac{ q_{12}}{-a} = \frac{ b_{12}b^{-1}}{-a}.
$$
Сравнивая диагональные коэффициенты при $\frac{1}{\lambda}$, получаем
$$
P_2A-AP_2+P^\prime_1 = QP_1,
$$
откуда следуют равенства  $(p_{11}^{1})\prime= q_{12}p_{21}^1$ и $(p_{22}^{1})^\prime = q_{21}p_{12}^1$, При должном выборе констант интегрирования получаем
\begin{equation*}
\begin{gathered}
p_{11}^1 = -\int_x^1 q_{12}(t)p_{21}^1(t)dt =  -\int_x^1 b_{12}b^{-1}(t)p_{21}^1(t)dt, \\
 p_{22}^1 = \int_0^x q_{21}(t)p_{12}^1(t)dt =  \int_0^x b_{21}b(t)p_{21}^1(t)dt
\end{gathered}
\end{equation*}
Тем самым, мы полностью определили матрицу $P_1(x)$. Заметим, что гладкость $P_1(x)$ совпадает с гладкостью коэффициентов первоначальной системы \eqref{MainMatrixEquation}, то есть $P_1(x) \in W_1^n[0, 1]$.

Теперь определим рекуррентную процедуру нахождения следующих $P_l(x)$. Предположим, что мы уже определили матрицу $P_{m-1}(x)$. Выпишем равенства для коэффициентов при $\lambda^{-m+1}$ и $\lambda^{-m}:$
\begin{equation}\label{m-1}
	P_m A-AP_m+P^\prime_{m-1} = QP_{m-1},
\end{equation}
\begin{equation}\label{m}
	P_{m+1} A-AP_{m+1}+P^\prime_{m} = QP_{m}.
\end{equation}
Записывая \eqref{m-1} в покоординатном виде и решая уравнения сначала для внедиагональных, а затем для диагональных элементов, получаем
\begin{equation}\label{p_21}
	p_{21}^{m} = \frac{-\frac{d}{dx}p_{21}^{m-1} + q_{21}p_{11}^{m-1}}{a} = \frac{-\frac{d}{dx}p_{21}^{m-1} - q_{21}\int_x^1 q_{12}(t)p_{21}^{m-1}(t)dt}{a},
\end{equation}
\begin{equation}\label{p_12}
	p_{12}^{m} = \frac{-\frac{d}{dx}p_{12}^{m-1} + q_{12}p_{22}^{m-1}}{-a} = \frac{-\frac{d}{dx}p_{12}^{m-1} + q_{12}\int_0^x q_{21}(t)p_{12}^{m-1}(t)dt}{-a}.
\end{equation}
Наконец, из уравнений для диагональных элементов \eqref{m} определяем $p_{11}^m$ и $p_{22}^m$:
\begin{equation}\label{p_11_and_p_22}
\begin{gathered}
	p_{11}^{m} = -\int_x^1 q_{12}(t)p_{21}^m(t)dt =  -\int_x^1 b_{12}b^{-1}(t)p_{21}^m(t)dt, \\
 p_{22}^{m} = \int_0^x q_{21}(t)p_{12}^m(t)dtt =  \int_0^x b_{21}b(t)p_{21}^1(t)dt.
\end{gathered}
\end{equation}
Из формул \eqref{p_21}-\eqref{p_11_and_p_22} следует, что гладкость функций в матрице $P_{m+1}$ на единицу ниже гладкости функций в матрице $P_{m}$. Поэтому корректное построение матрицы $P_{m+1}$ возможно только в том случае, если $P_{m}$ допускает взятие еще одной производной. Чтобы гарантировать абсолютную непрерывность итоговой матрицы, закончим ее построение на $m=n$-ом  шаге, так как $P_{n}(x)$  принадлежит пространству $W_1^{n-(n-1)}[0, 1] = W_1^1[0, 1]$.

\bigskip

\subsection{ Оценки остатка}
Для оценки остатка в асимптотическом представлении $Y(x, \lambda)$ определим набор из четырех интегралов
\begin{equation}\label{Vij}
\begin{gathered}
v_{11}(s, x, \lambda) = -\int_{\max\{x, s\}}^1 q_{12}(t)e^{-\lambda[\rho(t)-\rho(s)]}dt, \,\ v_{12}(s, x, \lambda) = -\int_{\max\{x, s\}}^1 q_{12}(t)e^{\lambda[\rho(x)-\rho(t)]}dt,
\\
v_{21}(s, x, \lambda) = -\int_0^{\min\{x, s\}} q_{21}(t)e^{-\lambda[\rho(x)-\rho(t)]}dt,  \,\
v_{22}(s, x, \lambda) = -\int_{0}^{\min\{x, s\}} q_{21}(t)e^{\lambda[\rho(t)-\rho(s)]}dt,
\end{gathered}
\end{equation}
где $q_{ij}$ определены в \eqref{Q}. Отметим, что в \eqref{Vij} функции $q_{12}, q_{21}$ принадлежат по меньшей мере  классу $L_1[0, 1]$, а с учетом пределов интегрирования действительная часть показателей экспонент ограничена сверху числом $\kappa\rho(1)$, так как $\Re\lambda \ge \kappa $ при $\lambda \in D$. В случае, когда коэффициенты $b_{21}(x), b_{12}(x) \in W_1^1[0, 1]$, а коэффициенты $b_{11}(x), b_{22}(x) \in L_1[0, 1]$, эти интегралы допускают интегрирование по частям
\begin{equation}\label{PartVij}
\begin{gathered}
v_{11}(s, x, \lambda) = \left.\lambda^{-1}\frac{q_{12}(t)}{a(t)}e^{-\lambda[\rho(t)-\rho(s)]}\right\vert_{t=\max\{x, s\}}^{t=1}+\lambda^{-1}\hat{v}_{11}(s, x, \lambda), \\
v_{12}(s, x, \lambda) = \left.\lambda^{-1}\frac{q_{12}(t)}{a(t)}e^{\lambda[\rho(x)-\rho(t)]}\right\vert_{t=\max\{x, s\}}^{t=1}+\lambda^{-1}\hat{v}_{12}(s, x, \lambda), \\
v_{21}(s, x, \lambda) = \left.-\lambda^{-1}\frac{q_{21}(t)}{a(t)}e^{-\lambda[\rho(x)-\rho(t)]}\right\vert_{t=0}^{t=\min\{x, s\}}+\lambda^{-1}\hat{v}_{21}(s, x, \lambda), \\
v_{22}(s, x, \lambda) = \left.-\lambda^{-1}\frac{q_{21}(t)}{a(t)}e^{\lambda[\rho(t)-\rho(s)]}\right\vert_{t=0}^{t=\min\{x, s\}}+\lambda^{-1}\hat{v}_{22}(s, x, \lambda),
\end{gathered}
\end{equation}
где
\begin{equation}\label{HatVij}
\begin{gathered}
\hat{v}_{11}(s, x, \lambda) = -\int_{\max\{x, s\}}^1 \left(\frac{q_{12}(t)}{a(t)}\right)' e^{-\lambda[\rho(t)-\rho(s)]}dt, \\ \hat{v}_{12}(s, x, \lambda) = -\int_{\max\{x, s\}}^1 \left(\frac{q_{12}(t)}{a(t)}\right)' e^{\lambda[\rho(x)-\rho(t)]}dt, \\
\hat{v}_{21}(s, x, \lambda) = \int_0^{\min\{x, s\}} \left(\frac{q_{21}(t)}{a(t)}\right)' e^{-\lambda[\rho(x)-\rho(t)]}dt, \\
\hat{v}_{22}(s, x, \lambda) = \int_0^{\min\{x, s\}} \left(\frac{q_{21}(t)}{a(t)}\right)' e^{\lambda[\rho(t)-\rho(s)]}dt.
\end{gathered}
\end{equation}

В случае интегрируемых коэффициентов для оценки остатков удобно использовать обозначения
\begin{equation}\label{Upsilon}
\Upsilon(\lambda) = \max\limits_{i, j, s, x} |v_{ij}(s, x, \lambda)|, \,\ C_q^{int} = \int_0^1|q_{12}(t)|+|q_{21}(t)|dt,
\end{equation}
а в случае абсолютно непрерывных коэффициентов дополнительно введем обозначения
\begin{equation}\label{HatUpsilon}
\hat{\Upsilon}(\lambda) = \max\limits_{i, j, s, x} |\hat{v}_{ij}(s, x, \lambda)|, \,\
C_q = 4e^{2\kappa\rho(1)}  (\max\limits_{x}\left|\frac{q_{12}(x)}{a(x)}\right| + \max\limits_{x}\left|\frac{q_{21}(x)}{a(x)}\right|).
\end{equation}
\begin{lemma}
	$\Upsilon(\lambda) \rightarrow 0$ при $\Pi_{\kappa} \ni \lambda \rightarrow \infty$. Более того, если коэффициенты уравнения \eqref{MainVector} абсолютно непрерывные, то $\hat{\Upsilon}(\lambda) \rightarrow 0$ и $\Upsilon(\lambda)  \le C_q|\lambda|^{-1}$ при $\Pi_{\kappa} \ni \lambda \rightarrow \infty$.
\end{lemma}
\begin{proof}
	Приведем доказательство для $|v_{11}(s, x, \lambda)|$. Для остальных $v_{ij}$ это делается полностью аналогично. Так как $\rho(t)$ - монотонная функция, то можно сделать замену $\xi = \rho(t)$,  $\xi \in [0, \rho(1)]$ в $v_{11}(s, x, \lambda)$, после которой интеграл примет вид
	\begin{equation}\label{lemma_1}
	\int_{\rho(\max\{x, s\})}^{\rho(1)} f_{11}(\xi)e^{-\lambda[\xi - \rho(s)]}d\xi, \,\ f_{11}(\xi) = \frac{q_{12}(t(\xi))}{a(d(\xi))}.
	\end{equation}
	Заметим, что $f_{11}(\xi) \in L_1[0, \rho(1)]$, так как
	$\int_0^{\rho(1)} \frac{q_{12}(t(\xi))}{a(t(\xi))} d\xi = \int_0^1 q_{12}(t)dt$. Экспонента в \eqref{lemma_1} ограничена при $\Pi_{\kappa} \ni \lambda \rightarrow \infty$ и не превосходит $e^{\kappa\rho(1)}$, причем при $\xi \ne \rho(s)$ она убывает или осциллирует в $\Pi_{\kappa}$. Поэтому доказательство завершается так же, как в лемме Римана-Лебега. Для малого фиксированного $\ep>0$ подберем непрерывно дифференцируемую функцию $\hat{f}_{11}$ такую, что $\int_0^{\rho(1)} |f_{11}(\xi) - \hat{f}_{11}(\xi)|d\xi < \ep$, тогда
	$$
|v_{11}(s, x, \lambda)| \le \ep e^{\kappa\rho(1)} + 	\left|\int_{\rho(\max\{x, s\})}^{\rho(1)} \hat{f}_{11}(\xi)e^{-\lambda[\xi - \rho(s)]}\, d\xi\right|.
$$
	В результате после интегрирования по частям интеграл в последнем неравенстве можно будет оценить сверху  величиной $\ep e^{\kappa\rho(1)}$ при достаточно больших значениях $|\lambda|$, тогда $|v_{11}(s, x, \lambda)| \le 2\ep e^{\kappa\rho(1)}$, что завершает доказательство первого утверждения леммы.
	
	В случае $A(x), B(x) \in AC[0, 1]$ очевидно, что $q_{ij} \in AC[0,1]$, $\gamma_{ij} \in AC[0, 1]$ и $|a| > \varepsilon$ при некотором $\varepsilon >0$. Тогда $\frac{q_{ij}}{a} \in AC[0, 1]$ и  $\left(\frac{q_{ij}}{a}\right)' \in L_1[0, 1]$. Проводя аналогичные рассуждения получаем  $|\hat{v}_{ij}(s, x, \lambda)| \rightarrow 0$ при
	$\Pi_{\kappa} \ni \lambda \rightarrow \infty$ равномерно для всех $s, x \in[0,1]$.
	Из равенств \eqref{PartVij} видно, что $\Upsilon(\lambda) \le |\lambda^{-1}|(\frac{C_q}{2}+\hat{\Upsilon}(\lambda)) \le C_q|\lambda|^{-1}$, 
если $|\lambda|$  столь велико, что  выполнена оценка  $\hat{\Upsilon}(\lambda) \le \frac{C_q}{2}$.
\end{proof}

\bigskip

\subsection{Основная теорема}
Теперь мы можем сформулировать основную теорему работы, в которой докажем существование фундаментальной матрицы решений $Y(x, \lambda)$, имеющей представление \eqref{R} c явными формулам для матриц $R^m$.

Предварительно определим операторы, в терминах которых компактно перепишем формулы  \eqref{p_21}-\eqref{p_11_and_p_22}
\begin{gather*}
	(I_1 f)(x) = -\int_x^1 b_{12}(t) b^{-1}(t) f(t)\, dt,\quad (I_2 f)(x) = \int_0^x b_{21}(t) b(t) f(t)\, dt, \\
	(Df)(x) = \frac 1{a(x)} f'(x),  \ \  J_1 = \frac{b_{21}b}{a} I_1, \ \ J_2 = -\frac{b_{12}b^{-1}}a I_2,
\end{gather*}
где функции $a(x)$ и $b(x)$ определены в \eqref{ab}.
\begin{theorem*}
		Пусть  выполнено условие \eqref{a1} и все функции $a_{i}, b_{ij}$ принадлежат пространству $L_1[0, 1]$. Тогда при любом $\kappa \in \mathbb{R}$ существует фундаментальная матрица $Y(x, \lambda)$ уравнения \eqref{MainVectorEquation}, имеющая представление
		\begin{equation}
		Y(x, \lambda) = M(x)(I + R(x, \lambda)E(x, \lambda),
		\end{equation}
		где $M(x)$, $E(x, \lambda)$ определены в \eqref{ME}, a $R(x, \lambda)$ - голоморфная матриц-функция в полуплоскости $\Pi_{\kappa}^{+} = \{\lambda \in \mathbb{C} | \Re \lambda > -\kappa \}$ при достаточно больших $|\lambda|$, причем
для элементов этой матриц-функции  выполнены еценки
		$$\|r_{ij}(x, \lambda)\|_{C[0, 1]} \le C_0 \Upsilon(\lambda)$$
		с некоторой константой $C_0$.
		
	Если дополнительно выполнены условия \eqref{a3},
 то фундаментальную матрицу $Y(x, \lambda)$ можно выбрать такой, что $R(x, \lambda)$ допускает представление
	$$
	R(x, \lambda) = \frac{R^1(x)}{\lambda} + \dots + \frac{R^n(x)}{\lambda^n} + o(1)\lambda^{-n},
	$$
	где элементы матрицы $o(1)$ - бесконечно малые функции равномерно по $x \in [0, 1]$ при $\lambda \to \infty$, $\lambda \in \Pi_{\kappa}^{+}$.
	Матриц-функции $R^{m}$ вычисляются по формулам
	\begin{equation}\label{Formulas}
		\begin{aligned}
			&R^m = \begin{pmatrix}
				r_{11}^m & r_{12}^m \\[0.2cm]
				r_{21}^m & r_{22}^m
			\end{pmatrix}, \ \ \
		r_{11}^1 = I_1\frac{b_{21}b}{a}, \ \ \
 r_{21}^1 = \frac{b_{21}b}{a}, \ \ \
			r_{12}^1 = -\frac{b_{12}b^{-1}}{a}, \ \ \ 
r_{22}^1 = -I_2\frac{b_{12}b^{-1}}{a}, \\
			&r^{m+1}_{11} = I_1 r_{21}^{m+1}, \ \ \ r^{m+1}_{21} = (-D+ J_1)^m r^1_{21},  \ \ \
			r^{m+1}_{12} = (D+ J_2)^m r^1_{12},  \ \ \  r^{m+1}_{22} = I_2 r^{m+1}_{12}.
		\end{aligned}
	\end{equation}
Аналогичное утверждение верно, если $\Pi_{\kappa}^{+}$ заменить на   $\Pi_{\kappa}^{-} = \{\lambda \in \mathbb{C} | \Re \lambda < \kappa \}$.
\end{theorem*}

\bigskip

\subsection{Доказательство теоремы для случая $n=0$}

\begin{proof}
	Начнем с доказательства утверждения теоремы для интегрируемых коэффициентов. Вернемся к уравнению \eqref{MainMatrixEquationChange} для $Z(x, \lambda)$, возникающему после замены $Y(x, \lambda) \to Z(x, \lambda)$
$$
	Z^\prime = \lambda(AZ-ZA)+QZ.
$$
Перед формулировкой теоремы было построено формальное решение этого уравнения по формулам \eqref{Formulas}. Теперь мы предъявим решение уравнения \eqref{MainMatrixEquationChange}, существующее в $\Pi_{\kappa^+}^+$, и найдем его связь с формальным.

Запишем предыдущее уравнение в покомпонентном виде
	\begin{equation*}
	\begin{pmatrix}
	z_{11}^\prime & z_{12}^\prime \\[0.2cm]
	z_{21}^\prime & z_{22}^\prime
	\end{pmatrix}
	=
	\lambda\begin{pmatrix}
	0 & az_{12} \\[0.2cm]
	-az_{21} & 0
	\end{pmatrix}
	+
	\begin{pmatrix}
	q_{12}z_{21} & q_{12}z_{22} \\[0.2cm]
	q_{21}z_{11} & q_{21}z_{12}
	\end{pmatrix}.
	\end{equation*}
	Рассмотрим по отдельности первый и второй столбцы матрицы $Z(x, \lambda)$ и проинтегрируем с условиями $z_{1k}(1, \lambda) = \delta_{1k}$, $z_{2k}(0, \lambda) = \delta_{2k}$, где $\delta_{ij}$ - символ Кронекера. Тогда
	\begin{equation*}
	\begin{aligned}
	&\begin{pmatrix}
	z_{11} \\[0.2cm]
	z_{21}
	\end{pmatrix}
	-
	\begin{pmatrix}
	1 \\[0.2cm]
	0
	\end{pmatrix}
	=
	\begin{pmatrix}
	-\int_{x}^1 q_{12}(t) z_{21}(t, \lambda)dt \\[0.2cm]
	 \int_{0}^xq_{21}(t)e^{-\lambda[\rho(x)-\rho(t)]}z_{11}(t, \lambda)dt
	\end{pmatrix},
	\\
	&\begin{pmatrix}
	z_{12} \\[0.2cm]
	z_{22}
	\end{pmatrix}
	-
	\begin{pmatrix}
	0 \\[0.2cm]
	1
	\end{pmatrix}
	=
	\begin{pmatrix}
	-\int_x^1 q_{12}(t) e^{\lambda[\rho(x)-\rho(t)]}z_{22}(t, \lambda)dt \\[0.2cm]
	\int_{0}^x q_{21}(t) z_{12}(t, \lambda)dt
	\end{pmatrix}.
	\end{aligned}
	\end{equation*}
	Через $z_k$ обозначим $k$-ый столбец матрицы $Z(x, \lambda)$, а через $V_k(\lambda)$ интегральный оператор, определенный соответствующей правой частью предыдущих равенств. Тогда уравнение для $z_k$ запишется в виде
	\begin{equation}\label{z_k}
	z_k = z_k^0 + V_kz_k,
	\end{equation}
	где $z_k^0 = 	\begin{pmatrix}
	\delta_{1k} \\
	\delta_{2k}
	\end{pmatrix}$, $k =1, 2$.
	\begin{lemma}
		Оператор $V_k(\lambda): L_{\infty}\times L_{\infty} \rightarrow L_{\infty}\times L_{\infty}$ непрерывен, и $\|V_k(\lambda)\|_{L_{\infty} \rightarrow L_{\infty}} \le e^{2\kappa\rho(1)}C_q^{int}$. Более того, оператор $V_k^2(\lambda)$ является сжимающим, а именно $\|V_k^2(\lambda)\|_{L_{\infty} \rightarrow L_{\infty}} \le C_q^{int}\Upsilon(\lambda)$, где $\Upsilon(\lambda)$ определено в \eqref{Upsilon} и $\Upsilon(\lambda) \rightarrow 0$ при $D \ni \lambda \rightarrow \infty$ согласно Лемме 1. Кроме того, в случае абсолютно непрерывных коэффициентов оценку можно усилить $\|V_k^2(\lambda)\|_{L_{\infty} \rightarrow L_{\infty}} \le C_q^{int}\Upsilon  \le C |\lambda|^{-1}$.
	\end{lemma}
	\begin{proof}
		\begin{equation}\label{V1V2}
		V_1(\lambda)
		\begin{pmatrix}
		f_1 \\[0.2cm]
		f_2
		\end{pmatrix}
		=
		\begin{pmatrix}
		-\int_{x}^1 q_{12}(t) f_2(t)dt, \\[0.2cm]
		 \int_{0}^xq_{21}(t)e^{-\lambda[\rho(x)-\rho(t)]}f_1(t)dt
		\end{pmatrix}, \,\
		V_2(\lambda)
		\begin{pmatrix}
		f_1 \\[0.2cm]
		f_2
		\end{pmatrix}
		=
		\begin{pmatrix}
		- \int_x^1 q_{12}(t)e^{\lambda[\rho(x)-\rho(t)]}f_2(t)dt \\[0.2cm]
		\int_0^x q_{21}(t)f_1(t)dt
		\end{pmatrix}
		\end{equation}
		Из явного вида этих операторов следует, что 	 $\|V_k(\lambda)\|_{L_{\infty} \rightarrow L_{\infty}} \le e^{2\kappa\rho(1)}C_q^{int}$.
		\begin{equation}\label{V_1^2}
		V_1^2
		\begin{pmatrix}
		f_1 \\[0.2cm]
		f_2
		\end{pmatrix}
		=
		\begin{pmatrix}
		-\int_{x}^1 q_{12}(t)\int_0^t q_{21}(s)e^{-\lambda[\rho(t)-\rho(s)]}f_1(s)dsdt \\[0.2cm]
		-\int_{0}^x q_{21}(t)e^{-\lambda[\rho(x)-\rho(t)]}\int_t^1 q_{12}(s)f_2(s)dsdt
		\end{pmatrix}=
		\begin{pmatrix}
		\int_0^1 q_{21}(s)f_1(s)v_{11}(s, x, \lambda)ds \\[0.2cm]
		\int_0^1 q_{12}(s)f_2(s)v_{21}(s, x, \lambda)ds
		\end{pmatrix},
		\end{equation}
		\begin{equation}\label{V_2^2}
		V_2^2
		\begin{pmatrix}
		f_1 \\[0.2cm]
		f_2
		\end{pmatrix}
		=
		\begin{pmatrix}
		-\int_x^1 q_{12}(t)e^{\lambda[\rho(x)-\rho(t)]}\int_0^tq_{21}(s)f_1(s)dsdt \\[0.2cm]
		-\int_0^x q_{21}(t)\int_t^1(s)e^{\lambda[\rho(t)-\rho(s)]}f_2(s)dsdt
		\end{pmatrix}
		=
		\begin{pmatrix}
		\int_0^1 q_{21}(s)f_1(s)v_{12}(s, x, \lambda)ds \\[0.2cm]
		\int_0^1 q_{12}(s)f_2(s)v_{22}(s, x, \lambda)ds
		\end{pmatrix},
		\end{equation}
		где в \eqref{V_1^2}, \eqref{V_2^2} поменяли порядок интегрирования. Из последних равенств в \eqref{V_1^2}, \eqref{V_2^2} следует, что
		$$
		\|V_k^2
		\begin{pmatrix}
		f_1  \\f_2
		\end{pmatrix}\|_{L_{\infty} \times L_{\infty}} \le C_q^{int}\Upsilon(\lambda) 	 \left\|
		\begin{pmatrix}
		f_1  \\f_2
		\end{pmatrix}\right\|_{L_{\infty} \times L_{\infty}}.
		$$
	\end{proof}
	Вернемся к уравнению \eqref{z_k}. Представим его решение в виде формального ряда
	\begin{equation}\label{FormalSeries}
	z_k = z_k^0+ \left(\sum\limits_{\nu = 0}^\infty V_k^\nu(\lambda)\right) V_k(\lambda)z_k^0
	\end{equation}
	и воспользуемся результатами Леммы 2. Тогда при $\lambda$ достаточно больших имеем $\|V^2_k(\lambda)\| < \frac{1}{2}$, следовательно,
	\begin{multline*}
\|\sum\limits_{\nu = 0}^\infty V_k^\nu(\lambda)\|_{L_{\infty} \rightarrow L_{\infty}} =  \|\sum\limits_{\nu=0}^\infty V_k^{2\nu}(\lambda) \cdot [I+V_k(\lambda)]\| \\
 \le \sum\limits_{\nu=0}^\infty \|V_k^{2\nu}(\lambda)\| \cdot \|[I+V_k(\lambda)]\| <  2(1+e^{2\kappa\rho(1)}C_q^{int}).
 \end{multline*}
 Тем самым, получаем сходимость ряда \eqref{FormalSeries} по норме пространства $L_{\infty} \times L_{\infty}$. 
 
 Обозначим $C_V = 2(1+e^{2\kappa\rho(1)}C_q^{int})$. Заметим, что 
	\begin{equation*}
	V_1
	\begin{pmatrix}
	1 \\
	0
	\end{pmatrix}
	=
	\begin{pmatrix}
	0 \\
	-v_{21}(1, x, \lambda)
	\end{pmatrix},
	\,\
	V_2
	\begin{pmatrix}
	0 \\
	1
	\end{pmatrix}
	=
	\begin{pmatrix}
	v_{12}(0, x, \lambda) \\
	0
	\end{pmatrix}.
	\end{equation*}
	Из этих выражений следует, что $\|V_kz_k^0\|_{L_{\infty}\times L_{\infty}} \le \Upsilon(\lambda)$. Тогда из представления \eqref{FormalSeries} получаем
	$$
	\|z_k-z_k^0\|_{L_{\infty}\times L_{\infty}} \le \|[\sum\limits_{\nu = 0}^\infty V_k^\nu(\lambda)]\|_{L_{\infty} \rightarrow L_{\infty}} \cdot \|V_kz_k^0\|_{L_{\infty}\times L_{\infty}} \le C_V\Upsilon(\lambda),
	$$
 После перехода от матрицы $Z(x, \lambda)$ к матрице $Y(x, \lambda)$ получаем  утверждение теоремы для интегрируемых коэффициентов (для случая $n=0$).
\begin{note}
	Из вида матрицы $Y(x, \lambda)$ следует, что ее определитель имеет вид
	$$det(Y(x, \lambda)) = det(M(x))(1+o(1))det(E(x, \lambda)) \text{ при } \lambda \rightarrow \infty, \lambda \in \Pi_{\kappa}.$$ В частности, этот определитель отличен от нуля, и во всех рассматриваемых в работе случаях гладкости матрица $Y(x, \lambda)$ действительно является фундаментальной матрицей решений системы \eqref{MainVectorEquation} в области $\Pi_{\kappa}$.
\end{note}

\bigskip

\subsection{Доказательство теоремы для случая $n\geqslant 1$.}
Перейдем к доказательству второй части теоремы с коэффициентами подчиненными условиям \eqref{a2}. 
Запишем уравнение \eqref{MainMatrixEquationChange} в покомпонентном виде
	\begin{equation*}
	\begin{pmatrix}
	z_{11}^\prime & z_{12}^\prime \\
	z_{21}^\prime & z_{22}^\prime
	\end{pmatrix}
	=
	\lambda\begin{pmatrix}
	0 & az_{12} \\
	-az_{21} & 0
	\end{pmatrix}
	+
	\begin{pmatrix}
	q_{12}z_{21} & q_{12}z_{22} \\
	q_{21}z_{11} & q_{21}z_{12}
	\end{pmatrix}.
	\end{equation*}
	Ранее мы уже показали существование решений $z_k$ уравнения \eqref{MainMatrixEquationChange} вида
	$$
	z_k = \sum\limits_{\nu=0}^{\infty} V_k^{\nu}z_k^0,
	$$
	где $V_k$ определены в \eqref{V1V2}, а $z_k^0 = 	 \begin{pmatrix}
	\delta_{1k} \\
	\delta_{2k}
	\end{pmatrix}$, $k =1, 2$. Из следующего замечания следует, что для получения асимптотических формул с точностью $o(\lambda^{-n})$ нужно рассмотреть не весь ряд, а только его первые $2n+1$ слагаемых.
	\begin{note}
		Пусть матрицы $A, B \in AC[0, 1]$. Тогда
		$$
		z_k = \sum\limits_{\nu=0}^{\infty} V_k^{\nu}z_k^0 = \sum\limits_{\nu=0}^{2n} V_k^{\nu}z_k^0 + O_{kn}(x, \lambda),
		$$
		где $\|O_{kn}(x, \lambda)\|_{\L_{\infty}\times \L_{\infty}} = \|(\sum\limits_{\nu=0}^{\infty} V_k^{\nu}) \cdot V_k^{2n} \cdot V_kz_k^0\| \le C_V \Upsilon^{n+1}(\lambda) = O(\lambda^{-(n+1)})$.
	\end{note}
Таким образом,
	$$z_k = \sum\limits_{\nu=0}^{\infty} V_k^{\nu}z_k^0 =
	\sum\limits_{\nu=0}^{2n} V_k^{\nu}z_k^0 + o(\lambda^{-n}).
	$$	
	В случае $A, B \in W^n_1[0, 1]$ построенные нами формальные коэффициенты асимптотик при $\lambda^0, \lambda^1, \dots, \lambda^{n}$ содержатся в $\sum\limits_{\nu=0}^{2n}V_1^{\nu}	\begin{pmatrix}
	1 \\
	0
	\end{pmatrix}$ и $\sum\limits_{\nu=0}^{2n}V_2^{\nu}\begin{pmatrix}
	0 \\
	1
	\end{pmatrix}$. Однако, в этих частичных суммах также содержатся некоторые осциллирующие по $\lambda$ члены, возникающие из-за того, что $z_1$ и $z_2$ являются линейными комбинациями построенных формальных решений. Так как такие асимптотики плохо пригодны для изучения спектральных задач, дальше мы приведем процедуру получения из $z_1$ и $z_2$ асимптотических представлений, совпадающих с формальными.

	Сначала рассмотрим упрощенную ситуацию, когда  
	\begin{equation}\label{simple}
\begin{gathered}
	q_{21}(0) = q_{21}'(0) = \dots = q_{21}^{(n-1)}(0) = 0,
	\\
	q_{12}(1) = q_{12}'(1) = \dots = q_{12}^{(n-1)}(1) = 0.
\end{gathered}
\end{equation}
В этом случае осциллирующие члены не возникают. Это следует из того, что при дифференцировании по частям первые $n$ раз в операторах $V_1$ и $V_2$ не возникнет слагаемых, соответствующих значению подынтегральной функции в концевых точках (в точке $0$ для $V_1$, и в точке $1$ для $V_2$), а именно
	\begin{equation*}
\begin{aligned}
	V_1 \begin{pmatrix}
	1 \\
	0
	\end{pmatrix}
	&=
	\begin{pmatrix}
	0 \\
	\lambda^{-1}r_{21}^{1} + \lambda^{-2}(-D)r_{21}^{1}+ \dots + \lambda^{-n}(-D)^{n-1}r_{21}^{1} + o(\lambda^{-n})
	\end{pmatrix}
	\\
	V_2		
	\begin{pmatrix}
	0 \\
	1
	\end{pmatrix}
	&=
	\begin{pmatrix}
	\lambda^{-1}r_{12}^{1} + \lambda^{-2}Dr_{12}^{1} + \dots + \lambda^{-n}D^{n-1}r_{12}^{1} + o(\lambda^{-n})\\
	0
	\end{pmatrix}
	\\
V^2_1
	\begin{pmatrix}
	0 \\
	r_{21}
	\end{pmatrix}
	&=
	\begin{pmatrix}
	0 \\
	\lambda^{-1}J_1r_{21} +
	+\dots + \lambda^{-n} D^{n-1}J_1r_{21} - \lambda^{-n}e^{-\lambda \rho(x)} \int_0^x (D^{n-1}J_{1}p_{21}(t))'e^{\lambda \rho(t)}dt
	\end{pmatrix}
	\\ &=
	\begin{pmatrix}
	0 \\
	\lambda^{-1}J_{1}(r_{21}) + \lambda^{-2}DJ_{1}(r_{21})  +\dots + \lambda^{-n} D^{n-1}J_{1}r_{21} + o(\lambda^{-n})
	\end{pmatrix},
	\\ 
V^2_2
	\begin{pmatrix}
	r_{12} \\
	0
	\end{pmatrix}
	&=
	\begin{pmatrix}
	\lambda^{-1}J_{2}r_{12} + \dots + \lambda^{-n}D^{(n-1)}J_2r_{12} + \lambda^{-n}e^{\lambda \rho(x)} \int_x^{-1} (D^{n-1}J_{2}r_{12})'e^{-\lambda \rho(t)}dt\\
	0
	\end{pmatrix}
	\\ &=
	\begin{pmatrix}
	\lambda^{-1}J_{2}r_{12} + \dots + \lambda^{-n}D^{n-1}J_{2}p_{12} + o(\lambda^{-n})\\
	0
	\end{pmatrix},
	\end{aligned}
	\end{equation*}	
	где $(D^{n-1}J_{2}r_{12})'$, $((-D)^{n-1}J_{1}r_{21})' \in L_1[0, 1]$, а  
 оценки $o\left(\lambda^{-n}\right)$ следуют из леммы Римана-Лебега.

 Из последних четырех равенств следует, что в выражение $V_1^{2k-1}		
	\begin{pmatrix}
	1 \\
	0
	\end{pmatrix}$  войдут все слагаемые  из $(-D+J_1)^{k-1}(r_{21}^{1})$, получающиеся после раскрытия скобок, а также слагаемые с асимптотиками более высоких порядков. Аналогично, в выражение  $V_1^{2k}		
	\begin{pmatrix}
	1 \\
	0
	\end{pmatrix}$ войдут слагаемые из $I_{1}(-D+J_1)^{k-1}r_{21}^{1}$
и слагаемые более высоких порядков. Поэтому
\begin{equation*}
	\sum\limits_{\nu=0}^{2n}V_1^{\nu}
	\begin{pmatrix}
	1 \\
	0
	\end{pmatrix} = \sum\limits_{k=0}^n \lambda^{-k}
	\begin{pmatrix}
	r_{11}^{k} \\[0.2cm]
	r_{21}^{k}
	\end{pmatrix}
	+
	o(\lambda^{-n}), \,\
	\sum\limits_{\nu=0}^{2n}V_2^{\nu}
	\begin{pmatrix}
	0 \\
	1
	\end{pmatrix} = \sum\limits_{k=0}^n \lambda^{-k}
	\begin{pmatrix}
	r_{12}^{k} \\[0.2cm]
	r_{22}^{k}
	\end{pmatrix}
	+
	o(\lambda^{-n}).
\end{equation*}
	
В общем случае (без предположения \eqref{simple}) суммы $\sum\limits_{\nu=0}^{2n}V_1^{\nu}
	\begin{pmatrix}
	1 \\
	0
	\end{pmatrix}$ и $	 \sum\limits_{\nu=0}^{2n}V_2^{\nu}
	\begin{pmatrix}
	0 \\
	1
	\end{pmatrix}$ устроены сложнее. В этом случае имеем
	$$
	V_1
	\begin{pmatrix}
	1 \\
	0
	\end{pmatrix}
	=
	\begin{pmatrix}
	0 \\
	\lambda^{-1}r_{21}^{1} + \lambda^{-2}(-D)r_{21}^{1}+ \dots + \lambda^{-n}(-D)^{n-1}r_{21}^{1} + c_{1}(\lambda)e^{-\lambda \rho(x)} + o(\lambda^{-n})
	\end{pmatrix},
	$$
	где
	$
	c_1(\lambda) = -\sum\limits_{j=1}^n \lambda^{-j} (-D)^{j-1}(r_{21}^{1}) \rvert_{t=0},
	$
	$$
	V_2
	\begin{pmatrix}
	0 \\
	1
	\end{pmatrix}
	=
	\begin{pmatrix}
	\lambda^{-1}r_{12}^{1} + \lambda^{-2}Dr_{12}^{1} + \dots + \lambda^{-n}D^{n-1}r_{12}^{1} + c_2(\lambda)e^{\lambda \rho(x)} + o(\lambda^{-n}) \\
	0
	\end{pmatrix},
	$$
	где
	$
	c_2(\lambda) = -e^{\lambda \rho(1)}\sum\limits_{j=1}^n \lambda^{-j} D^{j-1}r_{12}^{1} \rvert_{t=1}.
	$
	Более того, такие дополнительные слагаемые возникают через каждые две степени оператора, а именно
	\begin{equation}\label{GeneralCaseZ1}
	\sum\limits_{\nu=0}^{2n}V_1^{\nu}
	\begin{pmatrix}
	1 \\
	0
	\end{pmatrix}
	=
	\sum\limits_{k=0}^n \lambda^{-k}
	\begin{pmatrix}
	r_{11}^{k} \\[0.2cm]
	r_{21}^{k}
	\end{pmatrix}
	+
	\sum\limits_{k=1}^n
	C_{1k}(\lambda)
	\sum\limits_{\nu=0}^{2n-(2k-1)}V_1^{\nu}
	\begin{pmatrix}
	0 \\
	e^{-\lambda \rho(x)}
	\end{pmatrix}
	+ o(\lambda^{-n}),
	\end{equation}
	\begin{equation}\label{GeneralCaseZ2}
	\sum\limits_{\nu=0}^{2n}V_2^{\nu}
	\begin{pmatrix}
	0 \\
	1
	\end{pmatrix}
	=
	\sum\limits_{k=0}^n \lambda^{-k}
	\begin{pmatrix}
	r_{12}^{k} \\[0.2cm]
	r_{22}^{k}
	\end{pmatrix}
	+
	\sum\limits_{k=1}^n
	C_{2k}(\lambda)
	\sum\limits_{\nu=0}^{2n-(2k-1)}V_2^{\nu}
	\begin{pmatrix}
	e^{\lambda \rho(x)} \\
	0
	\end{pmatrix}
	+ o(\lambda^{-n}),
	\end{equation}
	где
	$C_{1k}(\lambda) = \sum\limits_{j = k}^n \tilde{c}_{kj}\lambda^{-j}$, и 	 $C_{2k}(\lambda) = -e^{-\lambda \rho(1)}\sum\limits_{j=k}^n \bar{c}_{kj}\lambda^{-j}
	$.
	
	Бороться с этими осциллирующими слагаемыми нам поможет следующее 
важное замечание.
	\begin{note} Операторы $V_1$ и $V_2$ связаны следующими соотношениями
		$$
		V_1
		\left[e^{-\lambda \rho(x)}
		\begin{pmatrix}
		f_1(x, \lambda) \\[0.2cm]
		f_2(x, \lambda)
		\end{pmatrix}
		\right]
		=
		\begin{pmatrix}
		- e^{-\lambda \rho(x)} \int_x^1 q_{12}(t) e^{\lambda [\rho(x) - \rho(t)]}f_2(t, \lambda)dt \\[0.2cm]
		e^{-\lambda \rho(x)} \int_0^x q_{21}(t)f_1(t, \lambda)dt
		\end{pmatrix}
		=
		e^{-\lambda \rho(x)} V_2\begin{pmatrix}
		f_1(x, \lambda) \\[0.2cm]
		f_2(x, \lambda)
		\end{pmatrix},$$
		$$V_2
		\left[e^{\lambda \rho(x)}
		\begin{pmatrix}
		f_1(x, \lambda) \\[0.2cm]
		f_2(x, \lambda)
		\end{pmatrix}
		\right]
		=
		\begin{pmatrix}
		- e^{\lambda \rho(x)} \int_x^1 q_{12}(t) f_2(t, \lambda)dt \\[0.14cm]
		e^{\lambda \rho(x)} \int_0^x q_{21}(t)e^{-\lambda[\rho(x) - \rho(t)]}f_1(t, \lambda)dt
		\end{pmatrix}
		=
		e^{\lambda \rho(x)} V_1\begin{pmatrix}
		f_1(x, \lambda) \\[0.2cm]
		f_2(x, \lambda)
		\end{pmatrix}.$$
	\end{note}
Из замечания 3 следует, что 
\begin{equation*}
\begin{gathered}
\sum\limits_{\nu=0}^{2n-(2m-1)}V_1^{\nu}
	\begin{pmatrix}
	0 \\
	e^{-\lambda \rho(x)}
	\end{pmatrix}
	=
	e^{-\lambda \rho(x)}\sum\limits_{\nu=0}^{2n-(2m-1)}V_2^{\nu}
	\begin{pmatrix}
	0 \\
	1
	\end{pmatrix}, \\ 
\sum\limits_{\nu=0}^{2n-(2m-1)}V_2^{\nu}
	\begin{pmatrix}
	e^{\lambda \rho(x)} \\
	0
	\end{pmatrix}
	=
		e^{\lambda \rho(x)}\sum\limits_{\nu=0}^{2n-(2m-1)}V_1^{\nu}
	\begin{pmatrix}
	1 \\
	0
	\end{pmatrix}.
\end{gathered}
\end{equation*}
	Подставляя эти выражения в \eqref{GeneralCaseZ1}-\eqref{GeneralCaseZ2},
получаем
	\begin{equation}\label{GeneralCaseZ1Final}
	\sum\limits_{\nu=0}^{2n}V_1^{\nu}
	\begin{pmatrix}
	1 \\
	0
	\end{pmatrix}
	=
	\sum\limits_{k=0}^n \lambda^{-k}
	\begin{pmatrix}
	r_{11}^{k} \\[0.2cm]
	r_{21}^{k}
	\end{pmatrix}
	+
	\hat{C}_1(\lambda)
	e^{-\lambda \rho(x)}
	\sum\limits_{\nu=0}^{2n}V_2^{\nu}
	\begin{pmatrix}
	0 \\
	1
	\end{pmatrix}
	+
	o(\lambda^{-n}),
	\end{equation}
	\begin{equation}\label{GeneralCaseZ2Final}
	\sum\limits_{\nu=0}^{2n}V_2^{\nu}
	\begin{pmatrix}
	0 \\
	1
	\end{pmatrix}
	=
	\sum\limits_{k=0}^n \lambda^{-k}
	\begin{pmatrix}
	r_{12}^{k} \\[0.2cm]
	r_{22}^{k}
	\end{pmatrix}
	+
	\hat{C}_2(\lambda)
	e^{\lambda [\rho(x) - \rho(1)]}
	\sum\limits_{\nu=0}^{2n}V_1^{\nu}
	\begin{pmatrix}
	1 \\
	0
	\end{pmatrix}
	+
	o(\lambda^{-n}),
	\end{equation}
	где
	$\hat{C}_{1}(\lambda) = \sum\limits_{j = 1}^n \hat{c}_{1j}\lambda^{-j}$, $\hat{C}_{2}(\lambda) = \sum\limits_{j = 1}^n \hat{c}_{2j}\lambda^{-j}$ с некоторыми константами $\hat{c}_{1j}$, $\hat{c}_{1j}$,	а добавленные в последние слагаемые степени операторов в силу оценок на $\|V_1^2\|_{L_{\infty}\rightarrow L_{\infty}}$ и выражений для $C_{1k}(\lambda)$ дают вклад только в $o(\lambda^{-n})$. Из равенств \eqref{GeneralCaseZ1Final}-\eqref{GeneralCaseZ2Final} явно видно, что полученные решения являются линейными комбинациями формальных решений с коэффициентами $\hat{C}_1(\lambda)$, $\hat{C}_2(\lambda)$. Таким образом, мы доказали существование фундаментальных решений $\begin{pmatrix}
	y_{11} \\
	y_{21}
	\end{pmatrix}$, $\begin{pmatrix}
	y_{12} \\
	y_{22}
	\end{pmatrix}$ с асимптотиками
	\begin{equation*}
	\begin{aligned}
	\begin{pmatrix}
	y_{11} \\[0.2cm]
	y_{21}
	\end{pmatrix}
	=
	M(x)\left[\sum\limits_{k=0}^n \lambda^{-k} \begin{pmatrix}
	r_{11}^{k} \\[0.2cm]
	r_{21}^{k}
	\end{pmatrix}
	+
	o(\lambda^{-n})\right]\cdot e^{\lambda A_1(x)}
	+
	\hat{C}_1(\lambda)
	M(x)\left[\sum\limits_{\nu=0}^{2n}V_2^{\nu}
	\begin{pmatrix}
	0 \\[0.2cm]
	1
	\end{pmatrix}\right] \cdot e^{\lambda A_2(x)},
	\end{aligned}
	\end{equation*}
	\begin{equation*}
	\begin{aligned}
	\begin{pmatrix}
	y_{12} \\
	y_{22}
	\end{pmatrix}
	=
	M(x)\left[\sum\limits_{k=0}^n \lambda^{-k} \begin{pmatrix}
	r_{12}^{k} \\[0.2cm]
	r_{22}^{k}
	\end{pmatrix}
	+
	o(\lambda^{-n})\right]\cdot e^{\lambda A_2(x)}
	+
	\hat{C}_2(\lambda) e^{-\lambda\rho(1)}
	M(x)\left[\sum\limits_{\nu=0}^{2n}V_1^{\nu}
	\begin{pmatrix}
	1 \\
	0
	\end{pmatrix}\right] \cdot e^{\lambda A_1(x)}.
	\end{aligned}
	\end{equation*}
	Для доказательства теоремы остается взять линейные комбинации этих решений
	\begin{equation*}
	\begin{aligned}
	\begin{pmatrix}
	\hat{y}_{11} \\
	\hat{y}_{21}
	\end{pmatrix}
	&=
	\begin{pmatrix}
	y_{11} \\
	y_{21}
	\end{pmatrix}
	- \hat{C}_1(\lambda)
	\begin{pmatrix}
	y_{12} \\
	y_{22}
	\end{pmatrix}
	=
	M(x)\left[\sum\limits_{k=0}^n \lambda^{-k} \begin{pmatrix}
	r_{11}^{k} \\[0.2cm]
	r_{21}^{k}
	\end{pmatrix} + o(\lambda^{-n})\right]\cdot e^{\lambda A_1(x)},
	\end{aligned}
	\end{equation*}
	\begin{equation*}
	\begin{aligned}
	\begin{pmatrix}
	\hat{y}_{12} \\
	\hat{y}_{22}
	\end{pmatrix}
	&=
	\begin{pmatrix}
	y_{12} \\
	y_{22}
	\end{pmatrix}
	- \hat{C}_2(\lambda)e^{-\lambda \rho(1)}
	\begin{pmatrix}
	y_{11} \\
	y_{21}
	\end{pmatrix}
	=
	M(x)\left[\sum\limits_{k=0}^n \lambda^{-k} \begin{pmatrix}
	r_{12}^{k} \\[0.2cm]
	r_{22}^{k}
	\end{pmatrix} + o(\lambda^{-n})\right]\cdot e^{\lambda A_2(x)}
	\end{aligned}
	\end{equation*}
	Благодаря ограниченности экспонент $e^{\lambda[A_2(x) - A_1(x)]}$ и $e^{\lambda[A_1(x)-A_2(x)-\rho(1)]}$ при $\lambda \in \Pi_{\kappa}$ после взятия таких линейных комбинаций оценка остатка в виде $o(\lambda^{-n})$ в каждом из фундаментальных решений сохраняется, что и было использовано в последних равенствах.
\end{proof}

\end{document}